\documentclass[preprint]{elsarticle}

\usepackage[numbers]{natbib}
\usepackage{bm}
\usepackage{fullpage}
\usepackage{amssymb}
\usepackage{amsmath}
\usepackage{latexsym}
\usepackage{graphicx}
\usepackage{enumerate}
 \usepackage{tikz}
 \usetikzlibrary{topaths,calc}

\newtheorem{lemma}{Lemma}

\newtheorem{proposition}{Proposition}
\newtheorem{theorem}{Theorem}

\newproof{pf}{Proof}
\newproof{pot}{Proof of Theorem~\ref{thm:mainGeom}}

\def\bara{B\'ar\'any}

\def\cara{Carath\'eodory}

\def\R{\mathbb{R}}

\def\S{\mathbf{S}}

\def\conv{\operatorname{conv}}

\def\zero{{\bf 0}}

\begin{document}

\begin{frontmatter}

\author{Pauline Sarrabezolles}
\address{Universit\'e Paris Est, CERMICS (ENPC), 6-8 avenue Blaise Pascal, Cit\'e Descartes, F-77455 Marne-la-Vall\'ee, France}
\ead{pauline.sarrabezolles@cermics.enpc.fr}

\title{The colourful simplicial depth conjecture}
\begin{abstract}
Given $d+1$ sets of points, or colours, $\S_1,\ldots,\S_{d+1}$ in $\R^d$, a {\em colourful simplex} is a set $T\subseteq\bigcup_{i=1}^{d+1}\S_i$ such that $|T\cap \S_i|\leq 1$, for all $i\in\{1,\ldots,d+1\}$. The colourful \cara{} theorem states that, if $\zero$ is in the convex hull of each $\S_i$, then there exists a colourful simplex $T$ containing $\zero$ in its convex hull. Deza, Huang, Stephen, and Terlaky ({\em Colourful simplicial depth}, Discrete Comput. Geom., {\bf 35}, 597--604 (2006)) conjectured that, when $|\S_i|=d+1$ for all $i\in\{1,\ldots,d+1\}$, there are always at least $d^2+1$ colourful simplices containing $\zero$ in their convex hulls. We prove this conjecture via a combinatorial approach.
\end{abstract}


\end{frontmatter}
\maketitle

\section{Introduction}

A {\em colourful point configuration} is a collection of $d+1$ sets of points $\S_1,\ldots,\S_{d+1}$ in $\mathbb R^d$. A {\em colourful simplex} is a subset $T$ of $\bigcup_{i=1}^{d+1}\S_i$ such that $|T\cap \S_i|\leq 1$.  
The colourful Carath\'eodory theorem, proved by B\'ar\'any in 1982~\cite{Bar82}, states that, given a colourful point configuration $\S_1,\ldots,\S_{d+1}$ in $\mathbb R^d$ such that $\zero\in\bigcap_{i=1}^{d+1}\conv(\S_i)$, there exists a colourful simplex $T$ containing $\zero$ in its convex hull. In the same paper, B\'ar\'any uses this theorem  combined with Tverberg's theorem to give a bound on simplicial depth. His argument motivated the following question: how many colourful simplices, at least, contain $\zero$ in their convex hulls?

Let $\mu(d)$ denote the minimal number of colourful simplices containing $\zero$ in their convex hulls over all colourful point configurations $\S_1,\ldots,\S_{d+1}$ in $\mathbb R^d$ such that $\zero\in\conv(\S_i)$ and $|\S_i|=d+1$ for $i=1,\ldots,d+1$. The colourful Carath\'eodory theorem states that $\mu(d)\geq 1$. The quantity $\mu(d)$ has been investigated by~\citet*{DHST06}. They proved that $2d\leq\mu(d)\leq d^2+1$ and conjectured that $\mu(d)=d^2+1$. 
Later \citet*{BM06} proved that $\mu(d)\geq\max\left(3d,\left\lceil\frac{d(d+1)}5\right\rceil\right)$ for $d\geq 3$, \citet*{ST06} proved that $\mu(d)\geq\left\lfloor\frac{(d+2)^2}4\right\rfloor$, and \citet*{DSX11} showed that $\mu(d)\geq\left\lceil\frac{(d+1)^2}2\right\rceil$.  \citet*{DMS14} improved the bound to $\frac 12d^2+\frac 72d-8$ for $d\geq 4$. This latter result was obtained using a combinatorial generalization of the colourful point configurations suggested by \bara{} and known as {\em octahedral systems}, see~\cite{DSX11}.

We use this combinatorial approach to prove the conjecture.

\begin{theorem}\label{thm:mainGeom}
 The equality $\mu(d)=d^2+1$ holds for every integer $d\geq 1$.
\end{theorem}

The outline of the paper goes as follows. Section~\ref{sec:preliminary} is divided into two parts. First we define the octahedral systems and show their link with the colourful point configurations. Second, we introduce one of our main tools: the decomposition of an octahedral system over some elementary octahedral systems called umbrellas. Section~\ref{sec:proof} is devoted to the proof of Theorem~\ref{thm:mainGeom}. 

\section{Preliminaries}\label{sec:preliminary}

\subsection{Octahedral systems}\label{subsec:oct_sys}

Let $V_1,\ldots,V_n$ be $n$ pairwise disjoint finite sets, each of size at least $2$. An {\em octahedral system} is a set $\Omega\subseteq V_1\times\cdots\times V_n$ satisfying the {\em parity condition}: the cardinality of $\Omega\cap (X_1\times\cdots\times X_n)$ is even if $X_i\subseteq V_i$ and $|X_i|=2$ for all $i\in\{1,\ldots,n\}$. We use the terminology of hypergraphs to describe  an octahedral system: the sets $V_i$ are the {\em classes}, the elements in $V_i$ are the {\em vertices}, and the $n$-tuples in $V_1\times\cdots\times V_n$ are the {\em edges}. An edge whose $i$th component is a vertex $x\in V_i$ is {\em incident with the vertex $x$}, and conversely. A vertex $x$ incident with no edges is {\em isolated}. A class $V_i$ is {\em covered} if each vertex of $V_i$ is incident with at least one edge. Finally, the set of edges incident with $x$ is denoted by $\delta_\Omega(x)$ and the {\em degree of $x$}, denoted by $\deg_\Omega(x)$, refers to $|\delta_\Omega(x)|$.

\begin{lemma}\label{lem:nonEmpty}
 In every nonempty octahedral system, at least one class is covered.
\end{lemma}
\begin{pf}
Consider an octahedral system $\Omega\subseteq V_1\times\cdots\times V_n$. Suppose that no classes are covered. There is at least one isolated vertex $x_i$ in each $V_i$. Hence, if there were an edge $(y_1,\ldots,y_n)$ in $\Omega$, then the parity condition would not be satisfied for $X_i=\{x_i,y_i\}$.
\end{pf}

Given a colourful point configuration $\S_1,\ldots,\S_{d+1}$, the Octahedron Lemma~\cite{BM06,DHST06} states that, for any $\S'_1\subseteq \S_1,\ldots,\S'_{d+1}\subseteq \S_{d+1}$, with $|\S'_1|=\cdots=|\S'_{d+1}|=2$, the number of colourful simplices generated by $\bigcup_{i=1}^{d+1}\S'_i$ and containing $\zero$ in their convex hulls is even. The hypergraph over $V_1\times\cdots\times V_n$ where $V_i$ is identified with $\S_i$ and whose edges are identified with the colourful simplices containing $\zero$ in their convex hulls is therefore an octahedral system. Furthermore, a strengthening of the colourful \cara{} Theorem, given in~\cite{Bar82}, states that if $\zero\in\bigcap_{i=1}^{d+1}\conv(\S_i)$, then each point of the colourful point configuration is in some colourful simplices containing $\zero$ in their convex hulls. Hence, in an octahedral system $\Omega$ arising from such a colourful point configuration, each class $V_i$ is covered.

\subsection{Decompositions}

The following proposition, proved in~\cite{DMS14}, states that the set of all octahedral systems is stable under the ``symmetric difference'' operation.

\begin{proposition}\label{prop:stable}
 Let $\Omega$ and $\Omega'$ be two octahedral systems over the same vertex set. $\Omega\triangle\Omega'$ is an octahedral system.
\end{proposition}
\begin{pf}
Let $\Omega''=\Omega\triangle\Omega'$. As $\Omega''$ is a subset of $V_1\times\cdots\times V_n$, we simply check that the parity condition is satisfied. Consider $X_1\subseteq V_1,\ldots,X_n\subseteq V_n$ with $|X_i|=2$ for $i=1,\ldots,n$. We have
$$ |\Omega''\cap (X_1\times\cdots\times X_n)|=|\Omega\cap (X_1\times\cdots\times X_n)|+|\Omega'\cap (X_1\times\cdots\times X_n)|-2|\Omega\cap\Omega'\cap (X_1\times\cdots\times X_n)|.$$
All the terms of the sum are even, which allows to conclude.
\end{pf}

We now present a family of specific octahedral systems we call {\em umbrellas}. An umbrella $U$ is a set of the form $\{x^{(1)}\}\times\cdots\times\{x^{(i-1)}\}\times V_i\times\{x^{(i+1)}\}\times\cdots\times\{x^{(n)}\}$, with $x^{(j)}\in V_j$ for $j\neq i$. The class $V_i$ covered in $U$ is called its {\em colour}. $T=(x^{(1)},\ldots,x^{(i-1)},x^{(i+1)},\ldots,x^{(n)})$ is its {\em transversal}. An umbrella is clearly an octahedral system over $V_1\times\cdots\times V_n$ and we have the following proposition.

\begin{proposition}\label{prop:propCol}
 Two umbrellas of the same colour have an edge in common if and only if they are equal.
\end{proposition}
\begin{pf}
 An umbrella is entirely determined by its colour $V_i$ and its transversal $T$. Therefore, if two umbrellas of the same colour have an edge in common, they necessarily have the same transversal, which implies that they are equal.
\end{pf}

It was implicitly proved in Section 3 of~\cite{DMS14} that any octahedral system can be described as a symmetric difference of umbrellas. In this paper, we describe an octahedral system as a symmetric difference of other octahedral systems to bound its cardinality. We now focus on octahedral systems where the size of each class is equal to the number of classes.

 Consider a nonempty octahedral system $\Omega\subseteq V_1\times\cdots\times V_n$ with $|V_i|=n$ for all $i\in\{1,\ldots,n\}$. Denote by $i_1$ the smallest $i\in\{1,\ldots,n\}$ such that $V_{i}$ is covered in $\Omega$ and order the vertices $\{x_1,\ldots,x_n\}$ of $V_{i_1}$ by increasing degree: $\deg_\Omega(x_1)\leq\cdots\leq\deg_\Omega(x_n)$. We define $\mathcal U$ to be the set of umbrellas of colour $V_{i_1}$ containing an edge of $\Omega$ incident with $x_1$ and $W=\triangle_{U\in\mathcal U}U$. Let $\Omega_j$ be the set of all edges in $\Omega\triangle W$ incident with $x_j$. Formally, 
$$\mathcal U=\{U: U\mbox{ umbrella of colour }V_{i_1}\mbox{ and }U\cap\delta_\Omega(x_1)\neq\emptyset\}\mbox{ and }\Omega_j=\delta_{\Omega\triangle W}(x_j).$$
 Note that $|\mathcal U|=\deg_\Omega(x_1)$. In the remaining of the paper we refer to $(\mathcal U,\Omega_2,\ldots,\Omega_n)$ as a {\em suitable decomposition}.
\begin{lemma}\label{lem:decompo} 
Let $(\mathcal U,\Omega_2,\ldots,\Omega_n)$ be a suitable decomposition and $W=\triangle_{U\in\mathcal U}U$. We have
\begin{enumerate}
 \item[\textup{(i)}] $\Omega_j\cap\Omega_\ell=\emptyset$, for all $j\neq \ell$ (they have no edge in common),
\item[\textup{(ii)}] $\Omega=W\triangle \Omega_2\triangle\cdots\triangle\Omega_n$,
\item[\textup{(iii)}] $\Omega_j$ is an octahedral system, for all $j$,
\item[\textup{(iv)}] $\deg_{\Omega}(x_j)\geq\max(|\mathcal U|,|\Omega_j|-|\Omega_j\cap W|)$ for all $j$.
\item[\textup{(v)}] If $V_i$ is not covered in $\Omega$, then $V_i$ is neither covered in $\Omega\triangle W$ nor in any $\Omega_j$.
\end{enumerate}
\end{lemma}

The terminology suitable decomposition is due to point (ii) of Lemma~\ref{lem:decompo}.
\begin{pf}[Proof of Lemma~\ref{lem:decompo}]
We first prove (i). The $i_1$th component of any edge in $\Omega_j$ is $x_j$. Therefore, $\Omega_j$ and $\Omega_\ell$ have no edge in common if $j\neq\ell$.

We then prove (ii). There are exactly $\deg_\Omega(x_1)$ umbrellas of colour $V_{i_1}$ containing an edge of $\Omega$ incident with $x_1$. As $W$ is the symmetric difference of these umbrellas, $x_1$ is isolated in $\Omega\triangle W$. Thus, $\Omega_2,\ldots,\Omega_n$ form a partition of the edges in $\Omega\triangle W$ and $\Omega\triangle W=\Omega_2\triangle\cdots\triangle\Omega_n$. Taking the symmetric difference of this equality with $W$ we obtain $\Omega=W\triangle\Omega_2\triangle\cdots\triangle\Omega_n$.

We now prove (iii). By definition, the $\Omega_j$'s are subsets of $V_1\times\cdots\times V_n$. It remains to prove that they satisfy the parity condition. Consider $X_i\subseteq V_i$ with $|X_i|=2$ for $i=1,\ldots,n$. If $X_{i_1}$ does not contain $x_j$, there are no edges in $\Omega_j$ induced by $X_1\times\cdots\times X_n$. If $X_{i_1}$ contains $x_j$, the edges in $\Omega_j$ induced by $X_1\times\cdots\times X_n$ are the ones induced by $X_1\times\cdots\times X_{i_1-1}\times\{x_j\}\times X_{i_1+1}\times\cdots\times X_n$. As $x_1$ is isolated in $\Omega\triangle W$, those edges are exactly the edges in $\Omega\triangle W$ induced by $X_1\times\cdots\times X_{i_1-1}\times\{x_1,x_j\}\times X_{i_1+1}\times\cdots\times X_n$. According to Proposition~\ref{prop:stable}, $W$ is an octahedral system and $\Omega\triangle W$ as well, hence there is an even number of edges.

We prove (iv). We have $|\mathcal U|=\deg_\Omega(x_1)\leq\deg_\Omega(x_j)$ for all $j\in\{1,\ldots,n\}$. Furthermore, by definition of the symmetric difference, we have $(\Omega_2\triangle\cdots\triangle\Omega_n)\setminus W\subseteq \Omega$. This inclusion becomes $(\Omega_2\setminus W)\triangle\cdots\triangle(\Omega_n\setminus W)\subseteq \Omega$. As two $\Omega_\ell$'s share no edges, $\Omega_j\setminus W\subseteq\Omega$ and thus $\Omega_j\setminus W\subseteq\delta_\Omega(x_j)$ for all $j\in\{2,\ldots,n\}$. We obtain
$$|\Omega_j|-|\Omega_j\cap W|\leq\deg_\Omega(x_j).$$

Finally to prove (v) it suffices to prove that a class $V_i$ not covered in $\Omega$ remains not covered in $\Omega\triangle W$. Indeed, if a class is covered in an $\Omega_j$, it is also covered in $\Omega\triangle W$, as no two $\Omega_\ell$'s have an edge in common. Consider $V_i$ not covered in $\Omega$. There is a vertex $x\in V_i$ incident with no edges in $\Omega$. In particular, there are no edges in $\Omega$ incident with $x_1$ and $x$. Therefore, the umbrellas in $\mathcal U$, which are defined by the edges incident with $x_1$, contain no edges incident with $x$. Hence, $x$ is isolated in $W=\triangle_{U\in\mathcal U}U$ and in $\Omega$. Finally, $x$ remains isolated in $\Omega\triangle W$.
\end{pf}

Unlike the suitable decomposition of $\Omega$, which is a decomposition over general octahedral systems, the decomposition given in the following lemma is over umbrellas.

\begin{lemma}\label{lem:elDecompo}
Consider an octahedral system $\Omega\subseteq V_1\times\cdots\times V_n$ with $|V_i|=n$ for all $i\in\{1,\ldots,n\}$. There exists a set of umbrellas $\mathcal D$, such that $\Omega=\triangle_{U\in\mathcal D}U$ and such that the following implication holds:
\begin{center}
$V_i$ is the colour of some $U\in\mathcal D$ $\Longrightarrow$ $V_i$ is covered in $\Omega$.
\end{center}
\end{lemma}

\begin{pf}
 The proof works by induction on the number of covered classes in $\Omega$. If no classes are covered, then, according to Lemma~\ref{lem:nonEmpty}, $\Omega$ is empty. 

Suppose now that $k$ classes are covered, with $k\geq 1$, and consider a suitable decomposition $(\mathcal U,\Omega_2,\ldots,\Omega_n)$ of $\Omega$. Denote by $W$ the symmetric difference $W=\triangle_{U\in\mathcal U}U$. According to Proposition~\ref{prop:stable}, $W$ is an octahedral system, and so is $\Omega\triangle W$. There are stricly fewer covered classes in $\Omega\triangle W$ than in $\Omega$. Indeed, in $\Omega\triangle W$, the class $V_{i_1}$ is no longer covered, since $x_1$ is isolated, and according to (v) of Lemma~\ref{lem:decompo}, a class not covered in $\Omega$ remains not covered in $\Omega\triangle W$. By induction, there exists a set $\mathcal D'$ of umbrellas such that $\Omega\triangle W=\triangle_{U\in\mathcal D'}U$, and such that if there is an umbrella of colour $V_i$ in $\mathcal D'$, then $V_i$ is covered in $\Omega\triangle W$. As the umbrellas in $\mathcal D'$ are not of colour $V_{i_1}$, we have $\mathcal U\cap\mathcal D'=\emptyset$. Therefore, $\Omega=(\triangle_{U\in\mathcal U}U)\triangle(\triangle_{U\in\mathcal D'}U)$ and the set $\mathcal D=\mathcal U\cup\mathcal D'$ satisfies the statement of the lemma.  
\end{pf}

\section{Proof of the main result}\label{sec:proof}

The following theorem gives a general lower bound on the cardinality of an octahedral system. Our main theorem is a corollary of it. 

\begin{theorem}\label{thm:main}
Let $\Omega\subseteq{V_1\times\cdots\times V_n}$ be an octahedral system with $|V_1|=\cdots=|V_n|=n\geq 2$. If $k\geq 1$ classes among the $V_i$'s are covered, then $$|\Omega|\geq k(n-2)+2.$$
\end{theorem}

Before proving this theorem, we show how the main theorem can be deduced from it.

\begin{pot}
The inequality $\mu(d)\leq d^2+1$ is proved in~\cite{DHST06}. Let $\S_1,\ldots,\S_{d+1}$ be a colourful point configuration in $\mathbb{R}^d$. As explained in Section~\ref{subsec:oct_sys}, the set $\Omega\subseteq{V_1\times\cdots\times V_{d+1}}$, with $V_i=\S_i$ for $i=1,\ldots,d+1$ and whose edges correspond to the colourful simplices containing $\zero$ in their convex hulls, is an octahedral system. According to~\cite[Theorem 2.3.]{Bar82}, all the classes are covered in this octahedral system. Applying Theorem~\ref{thm:main} with $k=n=d+1$ gives the lower bound: $\mu(d)\geq d^2+1$.
\end{pot}

The remainder of the section is devoted to the proof of Theorem~\ref{thm:main}. The proof distinguishes two cases, corresponding to the following Propositions~\ref{prop:cas1} and~\ref{prop:cas2}. We first prove these propositions.

\begin{proposition}
 \label{prop:cas1}
Consider an octahedral system $\Omega\subseteq V_1\times\cdots\times V_n$ with $|V_i|=n$ for all $i\in\{1,\ldots,n\}$ and a class $V_i$ covered in $\Omega$. If $\Omega$ can be written as a symmetric difference of umbrellas, none of them being of colour $V_i$, then $|\Omega|\geq n^2$. 
\end{proposition}
\begin{pf}
 Let $\mathcal D$ be a set of umbrellas such that there are no umbrellas of colour $V_i$ in $\mathcal D$ and $\Omega=\triangle_{U\in\mathcal D}U$.  Denote by $y_1,\ldots,y_n$ the vertices of $V_i$, and by $\mathcal Q_j$ the set of umbrellas in $\mathcal D$ incident with $y_j$ for each $j\in\{1,\ldots,n\}$. As $\mathcal D$ does not contain any umbrellas of colour $V_i$, the umbrellas in $\mathcal Q_j$ all have transversals with $i$th component equal to $y_j$. Denote by $Q_j$ the symmetric difference of the umbrellas in $\mathcal Q_j$. We have that $Q_j$ is an octahedral system, according to Proposition~\ref{prop:stable}, and that $\delta_{\Omega}(y_j)=Q_j$, $Q_j\neq\emptyset$, and $Q_j\cap Q_\ell=\emptyset$ for all $j\neq \ell$. According to Lemma~\ref{lem:nonEmpty}, at least one class is covered in $Q_j$ and hence $|Q_j|\geq n$. Therefore, we have 
$$|\Omega|=\sum_{j=1}^n\deg_{\Omega}(y_j)= \sum_{j=1}^n|Q_j|\geq n^2$$
\end{pf}

\begin{proposition}\label{prop:cas2}
 Consider an octahedral system $\Omega\subseteq V_1\times\cdots\times V_n$ with $|V_i|=n$ for all $i\in\{1,\ldots,n\}$ and a suitable decomposition $(\mathcal U,\Omega_2,\ldots,\Omega_n)$ of $\Omega$. Consider $\mathcal O\subseteq\{\Omega_2,\ldots,\Omega_n\}$ such that for each $\Omega_j\in\mathcal O$ there is a class $V_i$ covered in $\Omega_j$ and in no other $\Omega_\ell\in\mathcal O$. Denote by $\mathcal P\subseteq\mathcal O$ the set of umbrellas in $\mathcal O$. We have 
$$|\Omega|\geq |\mathcal U|(n-|\mathcal O|)+\sum_{\Omega_j\in\mathcal O}|\Omega_j|-|\mathcal U|(|\mathcal O|-|\mathcal P|)-|\mathcal U|-|\mathcal P|+1.$$
\end{proposition}

\begin{pf}
Let $W=\triangle_{U\in\mathcal U}U$. The number of edges in $\Omega$ is equal to $\sum_{j=1}^n\deg_\Omega(x_j)$. We bound $\deg_\Omega(x_j)$ by $|\mathcal U|$ for $j=1$ and if $\Omega_j\notin\mathcal O$ and by $|\Omega_j|-|\Omega_j\cap W|$ otherwise, see (iv) in Lemma~\ref{lem:decompo}. We obtain $$|\Omega|\geq|\mathcal U|(n-|\mathcal O|)+\sum_{\Omega_j\in \mathcal O}\left(|\Omega_j|-|\Omega_j\cap W|\right).$$ 
We introduce a graph $G=(\mathcal V,\mathcal E)$ defined as follows. We use the terminology {\em nodes} and {\em links} for $G$ in order to avoid confusion with the vertices and edges of $\Omega$. The nodes in $\mathcal V$ are identified with the umbrellas in $\mathcal U$ and the $\Omega_j$'s in $\mathcal O$: $\mathcal V=\mathcal U\cup\mathcal O$. There is a link in $\mathcal E$ between two nodes if the corresponding octahedral systems have an edge in common. $G$ is bipartite: indeed, two umbrellas in $\mathcal U$ are of the same colour $V_{i_1}$ and, according to Proposition~\ref{prop:propCol}, they do not have an edge in common. According to Lemma~\ref{lem:decompo}, two $\Omega_j$'s do not have an edge in common either. 

For $\Omega_j$ in $\mathcal O$, we have $|\Omega_j\cap W|=\sum_{U\in\mathcal U}|\Omega_j\cap U|=\deg_G(\Omega_j)$, note that here the degree is counted in $G$. The fact that the umbrellas in $\mathcal U$ are disjoint proves the first equality. The second equality is deduced from the facts that $\Omega_j$ has at most one edge in common with each umbrella in $\mathcal U$, the one incident with $x_j$, and that $\Omega_j$ has no neighbours in $\mathcal O$. We obtain the following bound
\begin{eqnarray*}
 |\Omega|&\geq &|\mathcal U|(n-|\mathcal O|)+\sum_{\Omega_j\in\mathcal O}\left(|\Omega_j|-\deg_G(\Omega_j)\right)\\
&=& |\mathcal U|(n-|\mathcal O|)+\sum_{\Omega_j\in\mathcal O}|\Omega_j|-\deg_G(\mathcal O\setminus\mathcal P)-\deg_G(\mathcal P).
\end{eqnarray*}

Again, for the equality, we use the fact that $G$ is bipartite. The number of links in $\mathcal E$ incident with a node in $\mathcal O\setminus\mathcal P$ is at most $|\mathcal U|$. Hence, $\deg_G(\mathcal O\setminus\mathcal P)\leq |\mathcal U|(|\mathcal O|-|\mathcal P|)$. It remains to bound $\deg_G(\mathcal P)$. Note that if $U$ is an umbrella in $\mathcal P$, it is the only umbrella of its colour in $\mathcal P$, otherwise it would contradict the property of $\mathcal O$. We now prove that there are no cycles induced by $\mathcal P\cup\mathcal U$ in $G$.

Suppose there is such a cycle $\mathcal C$ and consider an umbrella $U$ of $\mathcal P$ in this cycle. Denote its colour by $V_i$ and its neigbours in $\mathcal C$ by $L$ and $R$. As $G$ is simple, $L$ and $R$ are distinct. $L$ and $R$ are both in $\mathcal U$, and hence are of colour $V_{i_1}$ and do not have an edge in common. Therefore $U\cap L$ and $U\cap R$ do not have an edge in common either, which implies that the $i$th component of the transversals of $L$ and $R$ are distinct. Note that two umbrellas adjacent in $\mathcal C$, both of colour distinct from $V_i$, have necessarily transversals with the same $i$th component. Hence there must be another umbrella of colour $V_i$ in the path in $\mathcal C$ between $L$ and $R$ not containing $U$. This is a contradiction since $U$ is the only umbrella in $\mathcal P$ of colour $V_i$.

 The number of links in $\mathcal E$ incident with $\mathcal P$ is then at most $|\mathcal U|+|\mathcal P|-1$. This allows us to conclude.
\end{pf}

\begin{pf}[Proof of Theorem~\ref{thm:main}]
Let $\Omega\subseteq{V_1\times\cdots\times V_n}$ be an octahedral system with $|V_1|=\cdots=|V_n|=n\geq 2$, and suppose that $k\geq 1$ classes $V_{i_1},\ldots,V_{i_k}$, with $i_1<\cdots<i_k$, are covered in $\Omega$. The proof works by induction on $k$.\\

If $k=1$, then $\Omega$ must contain at least $n$ edges for one class to be covered.\\

Assume now that $k>1$. If $|\mathcal U|\geq n-1$, then, according to (iv) of Lemma~\ref{lem:decompo}, $|\Omega|=\sum_{j=1}^n\deg_\Omega(x_j)\geq n|\mathcal U|\geq k(n-2)+2$ and we are done. Assume now that $|\mathcal U|\leq n-2$. We consider a suitable decomposition $(\mathcal U,\Omega_2,\ldots,\Omega_n)$ of $\Omega$ and distinguish two cases.

Case $1$: {\em One of the covered classes $V_i$, for $i\in\{i_2,\ldots,i_k\}$, is not covered in any $\Omega_j$}. Let $V_i$ be a covered class in $\Omega$, which is not covered in any $\Omega_j$. For each $j\in\{2,\ldots,n\}$, applying Lemma~\ref{lem:elDecompo} on $\Omega_j$ gives a set $\mathcal D_j$ of umbrellas, all of colour distinct from $V_i$, such that $\Omega_j=\triangle_{U\in\mathcal D_j}U$. We obtain $\Omega=(\triangle_{U\in\mathcal U}U)\triangle(\triangle_{j=2}^n\triangle_{U\in\mathcal D_j}U)$, according to (ii) of Lemma~\ref{lem:decompo}. Thus, we can apply Proposition~\ref{prop:cas1} which ensures that $$|\Omega|\geq n^2\geq k(n-2)+2.$$ 

Case $2$: {\em Each covered class $V_i$, for $i\in\{i_2,\ldots,i_k$\}, is covered in at least one of the $\Omega_j$.} Choose a set $\mathcal O\subseteq\{\Omega_2,\ldots,\Omega_n\}$, minimal for inclusion, such that each covered class $V_i$, for $i\in\{i_2,\ldots,i_k\}$, is covered in at least one of the $\Omega_j\in\mathcal O$. Such a set $\mathcal O$ satisfies the statement of Proposition~\ref{prop:cas2}. Applying this proposition, we obtain 
$$|\Omega|\geq |\mathcal U|(n-|\mathcal O|)+\sum_{\Omega_j\in\mathcal O}|\Omega_j|-|\mathcal U|(|\mathcal O|-|\mathcal P|)-|\mathcal U|-|\mathcal P|+1.$$

We now bound $\sum_{\Omega_j\in\mathcal O}|\Omega_j|$. Let $k_j$ be the number of classes covered in $\Omega_j$. By minimality of $\mathcal O$, there is at least one class covered in each $\Omega_j\in\mathcal O$, and according to (v) of Lemma~\ref{lem:decompo} we have $k_j<k$, hence $1\leq k_j< k$. By induction, the cardinality of $\Omega_j$ is at least $k_j(n-2)+2$. This lower bound is not good enough for the $\Omega_j\notin\mathcal P$ such that $k_j=1$. We denote by $\mathcal A$ those $\Omega_j$'s. We explain now how to improve the lower bound for $\Omega_j\in\mathcal A$. Only one class is covered in $\Omega_j$ and $\Omega_j\notin \mathcal P$. According to Lemma~\ref{lem:elDecompo}, $\Omega_j$ can be written as a symmetric difference of distinct umbrellas of the same colour. According to Proposition~\ref{prop:propCol}, these umbrellas are pairwise disjoint and $|\Omega_j|$ is equal to $n$ times the number of umbrellas in this decomposition. Since $\Omega_j$ is not an umbrella itself, otherwise $\Omega_j$ would have been in $\mathcal P$, there are at least two umbrellas in this decomposition. We obtain
$$\sum_{\Omega_j\in\mathcal O}|\Omega_j|\geq\left(\sum_{\Omega_j\in\mathcal O\setminus\mathcal A}k_j\right)(n-2)+2|\mathcal O\setminus\mathcal A|+2n|\mathcal A|=\left(\sum_{\Omega_j\in\mathcal O}k_j\right)(n-2)+2|\mathcal O|+n|\mathcal A|$$
We have thus $$|\Omega|\geq|\mathcal U|(n-|\mathcal O|)+\left(\sum_{\Omega_j\in\mathcal O}k_j\right)(n-2)+2|\mathcal O|+n|\mathcal A|-|\mathcal U|(|\mathcal O|-|\mathcal P|)-|\mathcal U|-|\mathcal P|+1.$$

\noindent Finally, we have
\begin{eqnarray}
2|\mathcal O|-|\mathcal P|-|\mathcal A| &\leq &\sum_{\Omega_j\in\mathcal O}k_j\label{eq:3}\\
k-1 &\leq & \sum_{\Omega_j\in\mathcal O}k_j\label{eq:2}
\end{eqnarray}
 Equation~\eqref{eq:3} is obtained by distinguishing the $\Omega_j$ with $k_j=1$ from those with $k_j\geq 2$. Equation~\eqref{eq:2} results from the fact that each class $V_{i_2},\ldots,V_{i_k}$ is covered in at least one $\Omega_j$ in $\mathcal O$. Thus,

\begin{eqnarray*}
|\Omega|&\geq& |\mathcal U|(n-|\mathcal O|)+\left(\sum_{\Omega_j\in\mathcal O}k_j\right)(n-2)+2|\mathcal O|+|\mathcal U||\mathcal A|-|\mathcal U|(|\mathcal O|-|\mathcal P|)-|\mathcal U|-|\mathcal P|+1\\
&\geq &(k-1)(n-2)+2|\mathcal O|-|\mathcal P|+1+\left(\sum_{\Omega_j\in\mathcal O}k_j-k+|\mathcal A|+n-2|\mathcal O|+|\mathcal P|\right)|\mathcal U|
\end{eqnarray*}

\noindent where we only used the inequalities $n\geq n-2\geq|\mathcal U|$ and~\eqref{eq:2}. According to~\eqref{eq:3}, the expression $$\left(\sum_{\Omega_j\in\mathcal O}k_j-k+|\mathcal A|+n-2|\mathcal O|+|\mathcal P|\right)$$ is nonnegative. Moreover, we have already noted that $|\mathcal U|=\deg_\Omega(x_1)$, which is at least $1$. Therefore,
$$|\Omega|\geq (k-1)(n-2)+2|\mathcal O|-|\mathcal P|+1+\sum_{\Omega_j\in\mathcal O}k_j-k+|\mathcal A|+n-2|\mathcal O|+|\mathcal P|.$$
Using~\eqref{eq:2} again, we obtain $$|\Omega|\geq k(n-2)+2.$$
\end{pf}

\section*{Aknowlegement}
The author thanks Antoine Deza for introducing her to the colourful simplicial depth conjecture and Fr\'ed\'eric Meunier for his thorough reading of the manuscript and his helpful comments.

\appendix


\begin{thebibliography}{00}
 \bibitem[I. \bara{}(1982)]{Bar82}
 I. B{\'a}r{\'a}ny, A generalization of {C}arath\'eodory's theorem, Disc. Math. 40 (1982) 141--152.

\bibitem[Deza et. al(2006)]{DHST06}
 A. Deza, S. Huang, T. Stephen, T. Terlaky, Colourful simplicial depth, Discrete Comput. Geom. 35 (2006) 597--604.

\bibitem[I. \bara{} and J. Matou{\v s}ek(2007)]{BM06}
I. B{\'a}r{\'a}ny, J. Matou{\v s}ek, Quadratically many colorful simplices, SIAM J. Discrete Math. 21 (2007) 191--198.

\bibitem[Stephen and Thomas(2008)]{ST06}
T. Stephen, H. Thomas, A quadratic lower bound for colourful simplicial depth, Journal of Combinatorial Optimization 16 (2008) 324--327. 

\bibitem[Deza et. al(2011)]{DSX11}
A. Deza, T. Stephen, F. Xie, More colourful simplices, Discrete Comp. Geom. 45 (2011) 272--278.

\bibitem[Deza et. al(2014)]{DMS14}
A. Deza, F. Meunier, P. Sarrabezolles, A combinatorial approach to colourful simplicial depth, SIAM J. Discrete Math. 28 (2014) 306--322.
\end{thebibliography}

\end{document}